\documentclass[A4paper,10pt,twocolumn]{article}

\usepackage{graphicx}
\usepackage{amsmath}
\usepackage{amssymb}
\usepackage{amsfonts}
\usepackage{epsfig}
\usepackage{color}
\usepackage[bookmarks=true,colorlinks=true,hyperindex=true,linkcolor=red,citecolor=blue]{hyperref} 
\newcommand{\q}[1]{\vert #1 \rangle}
\newcommand{\qd}[1]{\langle #1 \vert}

\newcommand{\daag}{^{\dagger}}
\newcommand{\minou}{\text{-}}
\newcommand{\ba}{\text{\bf{a}}}
\newcommand{\bn}{\text{\bf{N}}}
\newcommand{\ban}{\text{\bf{a}}_{\bar n}}
\newcommand{\bnn}{\text{\bf{N}}_{\bar n}}
\newcommand{\Pn}{\text{\bf{P}}_{\bar n}}
\newcommand{\bH}{\text{\bf{H}}}
\newcommand{\bid}{\text{\bf{I}}}

\newcommand{\nmax}{\bar n}
\newcommand{\Dh}{\mathcal{D}_\mathcal{H}}
\newcommand{\TDH}{T\mathcal{D}_\mathcal{H}}
\newcommand{\trace}{\mathrm{Tr}}

\newcommand{\rhon}{\bar{\rho}_{\bar{n}}}
\newcommand{\rhonss}{\bar{\rho}_{\bar{n},\infty}}

\title{\LARGE \bf Contraction and stability analysis of steady-states for open quantum systems described by Lindblad differential equations}

\author{Pierre Rouchon and Alain Sarlette\thanks{P.~Rouchon is with the Centre Automatique et Syst\`emes, Math\'{e}matiques et Syst\`{e}mes,  Mines ParisTech, 75272 Paris cedex 06, France. A.~Sarlette is with the SYSTeMS research group at Ghent University, 914 Technologiepark, 9052 Zwijnaarde, Belgium. {\tt\footnotesize pierre.rouchon@mines-paristech.fr, alain.sarlette@ugent.be}}
\thanks{This paper presents research results of the Belgian Network DYSCO (Dynamical Systems, Control, and Optimization), funded by the Interuniversity Attraction Poles Program, initiated by the Belgian State, Science Policy Office. The scientific responsibility rests with its authors. PR was partially supported by the ANR, Projet Blanc EMAQS ANR-2011-BS01-017-01. The authors want to thank Claude Le Bris for stimulating discussions on this subject.}%
}

\begin{document}

\maketitle

\begin{abstract}
For discrete-time systems, governed by Kraus maps, the work of D.~Petz has characterized the set of  universal contraction metrics. In the present paper, we use this characterization  to derive a set of quadratic Lyapunov functions for continuous-time systems, governed by Lindblad differential equations, that have a steady-state with full rank. An extremity of this set is given by the Bures metric, for which the quadratic Lyapunov function is obtained by inverting a Sylvester equation. We illustrate the method by providing a strict Lyapunov function for a Lindblad equation designed to stabilize a quantum electrodynamic ``cat'' state by reservoir engineering. In fact we prove that any Lindblad equation on the Hilbert space of the (truncated) harmonic oscillator, which has a full-rank equilibrium and which has, among its decoherence channels, a channel corresponding to the photon loss operator, globally converges to that equilibrium.
\end{abstract}


\section{Introduction}\label{sec:1}

The study of asymptotic convergence properties of open quantum systems has been reinvestigated in the last decades thanks to the development of quantum information~\cite{nielsen-chang-book} and reservoir engineering~\cite{PCZ96} to stabilize and protect fragile quantum states (see for example~\cite{KastoryanoThesis} and the references herein). The most standard models assume memoryless environments and are described by Lindblad differential equations (continuous-time) or Kraus maps (discrete-time). During the 70's, a line of work around \cite{Davies1970,Spohn1977,Frigerio1978} has developed a sufficient algebraic convergence criterion for these models, known as the Davies-Sphon-Frigerio (DSF) criterion. It ensures existence, uniqueness and global attractiveness of a full rank steady-state density operator. Meanwhile, the particular structure of the quantum master equations also allows to investigate convergence with a geometric approach. Indeed, one can define ``universal'' contraction metrics which are independent of the precise dynamics~\cite{PetzA}. As noticed in~\cite{LOHMIS1998A}, contraction analysis opens the way to developing systematic convergence-characterizing tools by examining where the contraction is strict.

\emph{In the present paper, we propose such a tool in the form of a systematic Lyapunov function design for quantum systems described by Lindblad differential equations; and we illustrate it by applying it to a practical reservoir engineering situation,} motivated by our recent reservoir engineering proposal in \cite{CatsA,CatsB}. We limit ourselves to finite dimension for simplicity, although in principle similar results should be possible in infinite dimension. We thereby start from a characterization of the metrics that are (non-strictly) contractive for all Kraus maps~\cite{PetzA}. This characterization then leads to a systematic construction of (non-strict) Lyapunov functions for all Lindblad equations with a full-rank equilibrium; a main advantage of these Lyapunov functions is that their time-derivative has an explicit expression which lends itself well to an analysis of critical points.  The resulting convergence criterion is different from the DSF criterion since it does not assume that the vector space spanned by the Lindblad operators is closed by Hermitian conjugation (modulo identity). This allows for instance to rapidly show the following result not covered by the DSF criterion, in the context of our engineered reservoir application: any Lindblad equation on the Hilbert space of the (truncated) harmonic oscillator, which has a full-rank equilibrium and a decoherence channel proportional to the photon loss operator, globally converges to that equilibrium.

The paper is organized as follows. Section \ref{sec:2} recalls the models describing the evolution of open quantum systems, both with discrete-time and continuous-time dissipative channels. Section \ref{sec:3} reviews known contractive metrics for discrete-time quantum systems, described by Kraus maps. Section \ref{sec:4} contains our main results about systematic Lyapunov functions for continuous-time quantum systems, described by Lindblad-Kossakowski differential equations. Section \ref{sec:5} applies the result to obtain a convergence proof for a reservoir engineered to stabilize quantum cat states of an electromagnetic field mode.\\


\section{Dynamics of open quantum systems}\label{sec:2}

All the material presented in this section can be found with much more details in \cite{HRbook} with developments centered around the experiment described in Section \ref{sec:5} and  in \cite{Alickibook,Tarasovbook,Breuerbook} for more formal theoretical exposures.

We consider a quantum system on a finite-dimensional Hilbert space $\mathcal{H}$, whose state is represented by a density operator $\rho$, i.e.~a self-adjoint, unit-trace positive semi-definite linear operator on $\mathcal{H}$ which represents quantum probabilities. We denote the set of all density operators on $\mathcal{H}$ by $\Dh$. If the quantum system is closed, it evolves according to the \emph{Hamiltonian} dynamics (Schr\"{o}dinger dynamics with $\hbar=1$)
$$\tfrac{d}{dt} \rho = -i [H,\rho]$$
where $i=\sqrt{-1}$, $[A,B] = AB-BA$ and the Hamiltonian $H$ is a self-adjoint operator on $\mathcal{H}$. Such evolution leads to a unitary propagator, i.e.~the evolution is an isometry for each $t$ and two states never converge towards each other.

This behavior can be changed by considering open quantum systems. In particular a measurement operation on a quantum system perturbs it in a non-unitary (in fact stochastic) way. The present paper in contrast considers open systems where information is lost to the environment without being measured. The latter can be conceived as the reduction to our target system, on $\mathcal{H}$, of a Hamiltonian dynamics on a larger Hilbert space $\hat{\mathcal{H}} = \mathcal{H} \otimes \mathcal{H^B}$ which is the tensor product of the target system and an external ``bath''. Open quantum systems can both be helpful, for stabilization, and detrimental if the bath induces dissipation ``in the wrong direction'' w.r.t.~our goal, also known as \emph{decoherence}. Typically, a quantum system interacts with both undesired decoherence sources, and stabilizing baths designed to counter them.

One of the mainstream open quantum models assumes that the environment loses any information about the system between consecutive interactions. In discrete-time, this Markov assumption leads to a \emph{trace-preserving completely positive map} or \emph{Kraus map} $\Phi$ for the density operator of the target system \cite{Kraus83}:
\begin{eqnarray}\label{eq:Krausmap}
	\rho & \rightarrow & \Phi(\rho) = \sum_k \, M_k \, \rho \, M_k\daag \\
\label{eq:Krauscond} & & \text{where } {\textstyle \sum_k \;} M_k\daag M_k = \bid \, .
\end{eqnarray}
The $M_k$ are arbitrary linear operators on $\mathcal{H}$ (but satisfying condition \eqref{eq:Krauscond}), $\phantom{e}\daag$ denotes operator adjoint and $\bid$ is the identity operator. In continuous-time, the Markov assumption leads to a class of models known as \emph{Lindblad(-Kossakowski) differential equations} \cite{Lindblad,Kossakowski}:
\begin{eqnarray}\label{eq:Lindbladmap}
	\tfrac{d}{dt}\rho = \mathcal{L}(\rho) & = & -i [H,\rho]  \\
	\nonumber & & - \tfrac{1}{2}\sum_k\, (L_k\daag L_k \rho + \rho L_k\daag L_k -2 L_k \rho L_k\daag) \, .
\end{eqnarray}
The $L_k$ are arbitrary linear operators on $\mathcal{H}$. In the following sections, we review contraction-like tools to assess the convergence of \eqref{eq:Krausmap} and develop how they lead to Lyapunov convergence arguments for \eqref{eq:Lindbladmap}. The stationary points of such dissipative evolutions are also known as \emph{pointer states} in the physics literature \cite{Zurek03}.

A particular case of Markovian open quantum systems is obtained through \emph{reservoir engineering}. We start with a discrete reservoir.
In this case, at each step, a `new' auxiliary system (thus carrying no memory about past system states) is prepared in a predetermined initial state $\rho^B(0)$ and brought into controlled Hamiltonian interaction with the target system for a time $T$; the coupled systems thus evolve as
$$\rho \otimes \rho^B(0) \;\; \rightarrow \;\; U_T \,(\rho \otimes \rho^B(0))\, U_T\daag$$
where $U_T = U(T)$ is a unitary operator on $\mathcal{H} \otimes \mathcal{H^B}$, solution at time $T$ of the controlled Schr\"odinger equation
$$\tfrac{d}{dt} U(t) = - i \,H_{\text{controlled}}(t) \, U$$
with initial state $U(0)=I$. Target and auxiliary system are in general entangled.
After interaction the auxiliary system is discarded; the state of the target system is then described by taking the partial trace of $U_T \,(\rho \otimes \rho^B(0))\, U_T\daag$ over the auxiliary system Hilbert space, and the resulting effect is described by a Kraus map of type \eqref{eq:Krausmap}, where the $M_k$ depend on $\rho^B(0)$ and the function $H_{\text{controlled}}(t)$ over $[0,T]$. Now a new auxiliary system and interaction can be initiated, to iterate the (potentially time-varying) Kraus map. With proper tuning, a stabilizing effect can be obtained (see Section \ref{sec:5}), without using measurement-based feedback.

In the limit of a very large number of short interactions per time unit ($T=dt \ll 1$ and $U(T)$ close to $\bid$), a continuous-time model of type \eqref{eq:Lindbladmap} can be obtained. The latter is also often viewed as the continuous interaction of the target system with one large infinite-dimensional system (e.g. thermal field environment).\\


\section{Contractive metrics for Kraus maps}\label{sec:3}

The content of this section is also presented with much more details in \cite{KastoryanoThesis} and relies essentially on the key contributions~\cite{PetzB,PetzA}.

\noindent \textbf{Definition III.1:} \emph{A distance function $d: \; \rho_1,\rho_2 \,\rightarrow\, d(\rho_1,\rho_2) \in \mathbb{R}$ defines a contractive metric for a Kraus map $\Phi$ if and only if  $\,d(\Phi(\rho_1),\Phi(\rho_2)) \leq d(\rho_1,\rho_2)\,$ for all $\rho_1,\rho_2$.}\vspace{2mm}

Motivated by popular quantum results, we first consider `distance functions' in a wide sense: they may be asymmetric with respect to their two arguments, sometimes strictly defined only for full rank $\rho_1$ and/or $\rho_2$. The latter point may be annoying when the goal of an evolution is to reach a \emph{pure state}, i.e.~a rank one density operator; in presence of disturbances, a small residual proportional to $\bid$ usually solves the situation.

Thanks to the particular structure of Markovian open quantum systems, there exist several metrics which are contractive for \emph{any} Kraus map; see e.g.~\cite{KastoryanoThesis} for a summary:
\begin{itemize}
	\item The trace distance $d_{tr}(\rho_1,\rho_2) = \tfrac{1}{2}\, \trace(\sqrt{(\rho_1-\rho_2)^2})$.
	\item The Bures distance $d_B(\rho_1,\rho_2) = \sqrt{1-F(\rho_1,\rho_2)}$, where $F(\rho_1,\rho_2) = \trace(\sqrt{\sqrt{\rho_1}\rho_2\sqrt{\rho_1}})$ is the quantum \emph{fidelity}.
	\item The Chernoff distance $d_C(\rho_1,\rho_2) = \sqrt{1-Q(\rho_1,\rho_2)}$, where $Q(\rho_1,\rho_2) = \min_{0 \leq s \leq 1} \, \trace(\rho_1^s\,\rho_2^{1-s})\,$.
	\item The quantum relative entropy $d_S(\rho_1,\rho_2) = \sqrt{\trace(\rho_1(\log \rho_1-\log \rho_2 ))}$.
	\item The quantum $\chi^2$-divergence $d_{\chi^2}(\rho_1,\rho_2) = \sqrt{\trace(\,(\rho_1-\rho_2)\rho_2^{-1/2} (\rho_1-\rho_2)\rho_2^{-1/2}\,)}$.
	\item Hilbert's projective cone metric: $d_H(\rho_1,\rho_2) = \log(\,\lambda_{\max}(\rho_2^{-1/2} \rho_1 \rho_2^{-1/2}) \,/\,\lambda_{\min}(\rho_2^{-1/2} \rho_1 \rho_2^{-1/2}) \,)$ if $\rho_1$ and $\rho_2$ have the same support, and otherwise $d_H(\rho_1,\rho_2)=\infty$.
\end{itemize}
All these metrics have their advantages and drawbacks, and none of them is always \emph{strictly} contracting.
The Hilbert metric for instance is symmetric but non-Riemannian,
and also contractive in the dual (a.k.a. Heisenberg picture): defining
\begin{equation}
	\Phi^*(X) = \sum_k \, M_k\daag X M_k
\end{equation}
for each self-adjoint operator $X$ on $\mathcal{H}$, such that $\trace(\Phi*(X) \, \rho) = \trace(X \, \Phi(\rho))$, we have that $\Phi^*(\bid) = \bid$ and that $d_H$ is contractive for $\Phi^*$: this shows that any $\Phi^*$ tends to bring the eigenvalues of any operator closer to each other. Furthermore, the contraction ratio for $\Phi$ is given by $\text{tanh}(\Delta(\Phi)/4)$ where $\Delta(\Phi) = \max\, \{d_H(\rho_1,\rho_2) \,:\, \rho_1,\rho_2 > 0 \} \,$. See \cite{SSR-CDC2010,Reeb2011} for recent illustrations.\\

When considering deviations from a target solution, it can be attractive to compare $\rho$ and $\rho+\delta\rho$ with a Riemannian metric on the set of strictly positive Hermitian operators of trace one. This set is a sub-manifold of  the vector space of Hermitian operators. For each $\rho \in \Dh$ with full rank, $\delta\rho$ belonging to the tangent space at $\rho$, is any trace-less Hermitian operator. By linearity of $\Phi$, Definition III.1 readily yields:\vspace{2mm}

\noindent \textbf{Definition III.2:} \emph{A Riemannian distance metric $\Vert \cdot \Vert_{\cdot}\, : \rho, \delta\rho \in \Dh \times \TDH \,\rightarrow\, \Vert \delta\rho \Vert_{\rho} \in \mathbb{R}$ is contractive for $\Phi$ if and only if (with slight abuse of notation) $\Vert \Phi(\delta\rho) \Vert_{\Phi(\rho)} \leq \Vert \delta\rho \Vert_{\rho}$ for all traceless Hermitian $\delta\rho$ and all $\rho>0$.}\vspace{2mm}

The following proposition results from the remarkable result by Petz \cite{PetzA,PetzB}.\vspace{2mm}

\noindent \textbf{Proposition III.3:} \emph{Consider the Riemannian metrics that are contractive for all (finite-dimensional) Kraus maps. The set of these metrics can be parameterized via the formula}
\begin{equation}\label{eq:PetzA1}
	\Vert \delta\rho \Vert_{\rho}^2 = \int_0^1 \, \trace(\,\delta\rho (\delta\omega_s + \delta\omega_s\daag)/2 \,) \, m(s) \, ds
\end{equation}
\emph{where $\delta\omega_s$ is solution of the Sylvester equation}
\begin{equation}\label{eq:PetzA2}
	s \, \delta\omega_s \rho + \rho \delta\omega_s = \delta\rho
\end{equation}
\emph{and $m_s \, ds$ is a positive finite measure.}

\underline{Proof:} In \cite{PetzA} the metrics are parameterized by \emph{standard operator monotone (decreasing) functions}, that is functions $f$ on the positive cone satisfying $f(A) \geq f(B)$ if $A\leq B$ ($A$, $B$ Hermitian positive operators), with $x\, f(x) = f(1/x)$ and $f(1)=1$. Such functions can be parameterized by positive finite measures $m(s)$ for $s\in[0,1]$ such that
\begin{equation}\label{eq:fx}
	f(x) = \frac{1}{2}\, \int_0^1 \, \left(\frac{1}{sx+1} + \frac{1}{s+x}\right) \, m(s) \, ds \, .
\end{equation}
Introducing right and left multiplication superoperators $L_\rho \sigma = \rho \sigma$ and $R_\rho \sigma = \sigma \rho$, the metrics in \cite{PetzA} are written:
$$\Vert \delta\rho \Vert_{\rho}^2 = \trace (\delta\rho \, f(L_\rho R_\rho^{-1}) R_\rho^{-1} \delta\rho ) \; .$$
Plugging \eqref{eq:fx} into this yields the announced form.\hfill$\square$\vspace{2mm}

The two extreme metrics in this set are as follows.
\begin{itemize}
	\item A measure $m(s)$ concentrated on $1$, gives the Bures metric $\trace(\,\delta\rho\,\delta\omega\,) = 2\trace(\,\delta\omega^2\,\rho \,)$ with hermitian $\delta\omega$ satisfying $\rho \delta\omega+\delta\omega\rho = \delta\rho\, .$ Its geodesic distance is the Bures distance $d_B$.
	\item A measure $m(s)$ concentrated on $0$, gives the symmetric metric $\trace(\,\delta\rho\, \rho^{-1} \delta\rho\,)$.\\
\end{itemize}


\section{Contraction-based Lyapunov functions for Lindblad dynamics}\label{sec:4}

Our goal now is to make use of these Kraus map contractive metrics to study convergence of Lindblad-Kossakowski differential equations.
In particular, we want to improve the DFS criterion~\cite{Davies1970,Spohn1977,Frigerio1978} which gives a sufficient algebraic condition characterizing the omega limit set. We therefore look for a systematic and strong Lyapunov function design method.
For standard characterizations of the distance to a target state $\rho_{\text{target}}$, like e.g.~fidelity $F(\rho,\rho_{\text{target}}) = \trace(\sqrt{\sqrt{\rho}\rho_{\text{target}}\sqrt{\rho}})$, the time-derivative can be  hard to compute in an operator setting. It is therefore meaningful to look for Lyapunov function constructions with handy computational properties. Our main results build up as follows.\vspace{2mm}

\noindent \textbf{Theorem IV.1:} \emph{Consider a full-rank density operator $\rho$ subject to any Lindblad equation of the form \eqref{eq:Lindbladmap}. Consider any of the Kraus map contractive metrics parameterized as in Proposition III.3. Then such metric is also contractive under the Lindblad equation as, for $\rho$ and $\delta\rho$ following solutions of the Lindblad equation, we have for each fixed $s$:}
\begin{multline}\label{eq:dlyap}
	\tfrac{d}{dt} \trace(\delta\rho \delta\omega_s ) = \tfrac{d}{dt} \trace(\delta\rho \delta\omega_s\daag )
	\\
  - s \, \sum_k\, \trace([\delta\omega,L_k] \rho [\delta\omega,L_k]\daag)
	\\  - \sum_k\, \trace([\delta\omega\daag,L_k] \rho [\delta\omega\daag,L_k]\daag) \; \leq 0
\end{multline}
\emph{with $\delta\omega_s(\rho)$ solution of \eqref{eq:PetzA2} at each time.}

\underline{Proof:} It suffices to establish \eqref{eq:dlyap}, since then the contraction of the metric follows from integrating the negative function defined there over $m(s) ds$.

Consider the function $f_s(\rho,\delta\rho) = \trace(\delta\rho \delta\omega_s)$. We have, using \eqref{eq:PetzA2} and its adjoint,
\begin{multline}\label{eq:4proof}
	 \tfrac{d}{dt} f_s(\rho,\delta\rho) =\trace(\tfrac{d\delta\rho}{dt} \delta\omega_s + \delta\rho \tfrac{d\delta\omega_s}{dt}) \\
= \trace(\,\tfrac{d\delta\rho}{dt} \delta\omega_s + \delta\omega_s\daag \, (\rho \tfrac{d\delta\omega_s}{dt} + s \tfrac{d\delta\omega_s}{dt} \rho) \,) \\
= \trace(\,\tfrac{d\delta\rho}{dt} \delta\omega_s + \delta\omega_s\daag \, (\tfrac{d\delta\rho}{dt}-s\delta\omega_s \tfrac{d\rho}{dt}) - \tfrac{d\rho}{dt}\delta\omega_s \,)\\
= \trace\big(\, \tfrac{d\delta\rho}{dt} \, (\delta\omega_s + \delta\omega_s\daag)) -\tfrac{d\rho}{dt}\, (s \delta\omega_s\daag \delta\omega_s + \delta\omega_s \delta\omega_s\daag) \,\big) .
\end{multline}
Now the Lindblad equation \eqref{eq:Lindbladmap} can be plugged in for $\rho$ and $\delta\rho$. The term in $[H,\rho]$ yields a term proportional to\linebreak $\trace( B - B\daag)$ with
$$B = \delta\rho H \delta\omega_s + \delta\rho H \delta\omega_s\daag - (s \delta\omega_s\daag \delta\omega_s + \delta\omega_s \delta\omega_s\daag) \rho H \, .$$
Replacing $\delta\rho$ by $\delta\omega_s\daag \rho + s \rho \delta\omega_s\daag$ in the first term and by $s \delta\omega_s \rho + \rho \delta\omega_s$ in the second term, one easily sees that $\trace(B) = \trace(B\daag)$. So the Hamiltonian term has no effect on $\tfrac{d}{dt} f_s(\rho,\delta\rho)$.

A similar procedure can be followed for the terms in $L_k$, for each $k$. Namely, once \eqref{eq:Lindbladmap} is plugged into the last line of \eqref{eq:4proof}, replace $\delta\rho$ respectively by its expression as a function of $\delta\omega_s,\rho$ in the terms containing $\delta\omega_s\daag$ and by the adjoint expression (as a function of $\delta\omega_s\daag,\rho$) in the terms containing $\delta\omega_s$. Then several terms simplify, grouping the remaining ones precisely yields the right hand side of \eqref{eq:dlyap}. Since the content of the trace operator there is self-adjoint, the same result must hold for $f_s\daag = \trace(\delta\rho \delta\omega_s\daag)$.
\hfill $\square$\vspace{2mm}

\noindent \textbf{Corollary IV.2:} \emph{In particular for the Bures metric, we have}
$$\tfrac{d}{dt}\trace(\delta\rho \delta\omega) \, = \, -2 \sum_k\trace([\delta\omega,L_k]\, \rho \, [\delta\omega,L_k]\daag) \leq 0$$
\emph{with $\rho \delta\omega + \delta\omega \rho= \delta\rho$.}\vspace{2mm}

This result can be used to systematically build a Lyapunov function for any Lindblad differential equation, for which an equilibrium $\rho_\infty$ (and nothing more) is known. Namely, we take $\rho_\infty$ as basis for the tangent vector $(\rho-\rho_\infty)$ in the above results.\vspace{2mm}

\noindent \textbf{Corollary IV.3:} \emph{If the equation \eqref{eq:Lindbladmap} has a full-rank equilibrium $\rho_\infty$, then whatever the form of $H$ and of the $L_k$, the function}\vspace{2mm}
\begin{eqnarray}
	\label{eq:LyapFunc0} &&V_\text{Bures}(\rho)= \trace(\rho_\infty\, G_\rho^2) \\
	\label{eq:LyapFunc1} && \text{ with }\rho_\infty G_\rho + G_\rho \rho_\infty = \rho-\rho_\infty
\end{eqnarray}
\emph{is a (non-strict) Lyapunov function for \eqref{eq:Lindbladmap} and:}
\begin{equation}\label{eq:LyapFuncDer}
	\tfrac{d}{dt} V_\text{Bures}(\rho) = - \sum_k \trace(\, [G_\rho,L_k]\, \rho_\infty \,[G_\rho,L_k]\daag \,) \, .
\end{equation}

In order to conclude about global convergence of the state to $\rho_\infty$, it then remains to examine how the commutator of the set of all $L_k$ relates to $\rho_\infty$, if necessary with a LaSalle-type argument.\\


\section{Application: a reservoir for quantum cats}\label{sec:5}

\subsection{Physical description}

\begin{figure}
	\includegraphics[width=80mm]{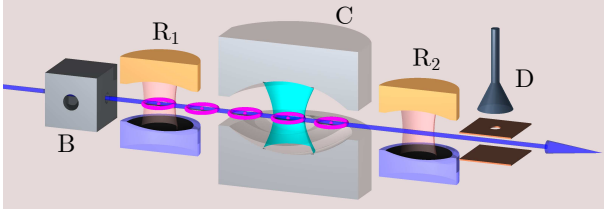}
	\caption{Simplified scheme of the cavity quantum electrodynamics experiment at Ecole Normale Sup\'erieure, Paris.}\label{fig:setup}
\end{figure}

We apply the above framework to a cavity quantum electrodynamics experiment, whose goal is to manipulate the (quantum) state of an electromagnetic field mode through its interaction with atoms. A scheme of the setup is shown on Figure \ref{fig:setup}; see \cite{HRbook} for a thorough explanation of its working.

In summary, a microwave field at frequency $\omega$ is trapped in the cavity $C$ made of superconducting mirrors. Atoms are sent one after the other through $C$, with the goal that the field interacts with an atomic transition at frequency $\omega_0 \approx \omega$. To this end, each atom is initially prepared by $B$ in one of the two states concerned by the relevant transition, and a classical microwave pulse in $R_1$ allows to put it in any superposition of those two states, parameterized by $u$. In addition to choosing with $u$ the initial state of each atom, we can control the experiment by imposing during each interaction a tailored evolution of $\delta(t) = \omega_0(t)-\omega$, the detuning between field and atomic transition frequency, with a good time resolution (via a Stark effect on the atom). Varying $\delta$ around $0$ allows to combine resonant and dispersive physical effects, on the way to producing truly quantum phenomena (see below).

The elements $R_2$ and $D$ can be used to detect the atomic state e.g.~for feedback purposes \cite{PboxNature11}. In the present case they play no role, since we use the atomic interactions as an engineered reservoir to control the field. This means that we consider the Kraus map associated to the field evolution when tracing over the atomic states. In addition to this (supposedly) stabilizing action, we want to consider the effect of decoherence due to spontaneous photon loss of the field, characterizing how the engineered reservoir allows to counter it. We therefore begin by describing the situation without decoherence in discrete-time (Kraus map), then derive its associated Lindblad equation, and add the traditional decoherence channel to this Lindblad equation. From there, we can compute an equilibrium and analyze convergence using the result of Section \ref{sec:5}.

\subsection{Goal; Kraus and Lindblad models}

An electromagnetic field mode at frequency $\omega$ is modeled as a quantum harmonic oscillator. Its Hilbert space is spanned by the orthonormal basis of \emph{Fock states} $\q{n}$ with $n=0,1,2,... $, which are the eigenstates of its closed-system Hamiltonian:
$$\bH_C = \sum_{n=0}^{+\infty} \; \hbar \omega (n+1/2) \, \q{n}\qd{n} = \bid/2 + \hbar\omega \bn \, .$$
Here we have defined the photon-number operator $\bn = \sum_{n=0}^{+\infty} n \q{n}\qd{n}$ and the identity operator $\bid = \sum_{n=0}^{+\infty} \q{n}\qd{n}$. We will denote the generic field state by $\rho$ or $\q{\psi}$.

The closest in this setting to a classical field state of complex amplitude $\alpha = A e^{i\phi}$, is the so-called \emph{coherent state}
$$\q{\alpha} = e^{-\vert \alpha \vert^2/2} \, \sum_{n=0}^{+\infty} \, \frac{\alpha^n}{\sqrt{n!}} \, \q{n} \, .$$
To highlight specifically quantum features, in \cite{CatsA,CatsB} our goal is instead to reach quantum superpositions of such coherent states, of the type
\begin{equation}\label{def:cat}
\q{\psi} = \sum_{k=1}^N  \beta_k \, \q{\alpha e^{k 2i\pi/N}}
\end{equation}
with fixed and known $\beta_k \in \mathbb{C}$, $\vert \beta_k \vert^2 = 1/N$. These are informally referred to as \emph{cat states} in honor of the famous Schr\"odinger cat. A particular set of cat states is generated by passing a coherent state through a \emph{Kerr medium}, i.e.~applying the Hamiltonian evolution corresponding to
$$\bH_K = \zeta_K \bn + \gamma_K \bn^2 \, .$$
Namely, $\q{\psi} = e^{-i t_K \bH_K}$ is a cat state of type \eqref{def:cat} with $N$ components if $t_K\gamma_K=\pi/N$. Unlike coherent states, cat states are very quickly destroyed by typical perturbations. Therefore, the in practice slow procedure of generating a cat state through a Kerr medium is poorly usable, and more subtle schemes must be devised to \emph{stabilize} them.

In \cite{CatsA,CatsB}, by tailoring the initial atom state ($u$) and imposing a simple time-varing profile of the interaction parameter ($\delta(t)$), we were able to propose an engineered reservoir such that one atom-field interaction, lasting a time $T$, applies the following Kraus map to the field (modulo some reasonable approximations; see \cite{CatsB} for details):
\begin{eqnarray}\label{eq:CatKraus}
	\rho & \rightarrow & M_1 \rho M_1\daag + M_2 \rho M_2\daag \quad \text{ with} \\
	\nonumber M_k & = & e^{-i h_\bn}\, \bar{M}_k \, e^{i h_\bn}\;\;\; \text{ for } k=1,2 \,;\\
	\nonumber \bar{M}_1 & = & \cos(\tfrac{u}{2}) \, \cos(\tfrac{\theta \sqrt{\bn}}{2}) + \sin(\tfrac{u}{2}) \, \frac{\sin (\tfrac{\theta \sqrt{\bn}}{2})}{\sqrt{\bn}} \, \ba\daag \,; \\
	\nonumber \bar{M}_2 & = & \sin(\tfrac{u}{2}) \, \cos(\tfrac{\theta \sqrt{\bn+\bid}}{2}) - \cos(\tfrac{u}{2}) \, \ba \, \frac{\sin (\tfrac{\theta \sqrt{\bn}}{2})}{\sqrt{\bn}} \,;\\
	\nonumber h_\bn & = & \phi\, \bn^2 + f(\phi)\, \bn \, .
\end{eqnarray}
Here the \emph{photon annihilation operator} is defined as $\ba = \sum_{n=0}^{+\infty} \, \sqrt{n} \q{n\minou 1}\qd{n}$ and satisfies $\ba\daag \ba = \bn$. Parameter $u \in [0,2\pi]$ reflects initial atomic state choice, while $\theta$ and $\phi$ can be varied in $[0,2\pi]$ by properly tailoring $\delta(t)$ during the interaction. The trick for ``simulating'' the Kerr-Hamiltonian-like factor $h_\bn$ is to play with the non-commutation of evolutions associated to the different values that $\delta(t)$ takes over time.\vspace{2mm}

\noindent \textbf{Lemma V.1:} \emph{The dynamics of $\rho$ under the Kraus map associated to $M_1,M_2$, corresponds to dynamics for $\bar{\rho} = e^{i h_\bn} \rho e^{-i h_\bn}$ according to the Kraus map associated to $\bar{M}_1,\bar{M}_2$. That is, the state trajectory $\rho(t)$ is equivalent to the state trajectory $\bar{\rho}(t)$ viewed through a Kerr medium.} \hfill $\square$\vspace{2mm}

Thus if $\bar{M}_1,\bar{M}_2$ stabilize a state $\bar{\rho}_\infty$ close to a coherent one -- this is a not too exotic task -- then the actual interaction can stabilize a state $\rho_\infty$ close to a cat state. In the following, for simplicity we focus on the case $\phi=\pi$ corresponding to a 2-component cat; a similar analysis can be carried out for a more general case.\\

If $u,\theta$ are small, the Kraus maps for $\rho$ and $\bar{\rho}$ are close to the identity and can be viewed as discretizations of Lindblad differential equations. One readily computes:
\begin{equation}\label{eq:cont-approx}
	\tfrac{d}{dt}\bar{\rho} = [\beta \ba\daag-\beta\daag \ba,\, \bar{\rho}] - \tfrac{\kappa}{2} (\bn \bar{\rho} + \bar{\rho} \bn - 2 \ba \bar{\rho} \ba\daag) \, ,
\end{equation}
with $\beta \, dt = u \theta / 4$ and $\kappa\, dt = \theta^2 / 4$. To this we can add Lindblad terms corresponding to disturbances. The dominant disturbance is spontaneous photon loss, corresponding to a term $- \tfrac{\kappa_c}{2} (\bn \rho + \rho \bn - 2 \ba \rho \ba\daag )$ in the original frame. Transforming that term to the frame of $\bar{\rho}$ and adding it to \eqref{eq:cont-approx} yields the dynamics:
\begin{eqnarray}
\nonumber	\tfrac{d}{dt}\bar{\rho} & = & \beta [ \ba\daag-\ba,\, \bar{\rho}]  - \tfrac{\kappa}{2} (\bn \bar{\rho} + \bar{\rho} \bn - 2 \ba \bar{\rho} \ba\daag )\\
 \label{eq:rhoKDecoherence} & &   - \tfrac{\kappa_c}{2} (\bn \bar{\rho} + \bar{\rho} \bn - 2 e^{i\pi\bn}\ba \bar{\rho} \ba\daag e^{-i\pi\bn} ) \, .\\
\nonumber & = & \beta [ \ba\daag-\ba,\, \bar{\rho}] - \tfrac{\kappa+\kappa_c}{2} (\bn \bar{\rho} + \bar{\rho} \bn - 2 \ba \bar{\rho} \ba\daag ) \\
\nonumber & & - \kappa_c (\ba \bar{\rho} \ba\daag - e^{i\pi\bn}\ba \bar{\rho} \ba\daag e^{-i\pi\bn} ) \; .
\end{eqnarray}
The last term induces a structural difference between equations \eqref{eq:cont-approx} and \eqref{eq:rhoKDecoherence}; due to the change of variables, it implies non-local interference since an operator is compared to its rotation by $\pi$ in electromagnetic phase space.

\subsection{Convergence analysis}

We analyze the convergence properties of $\bar{\rho}$ under \eqref{eq:rhoKDecoherence}; the corresponding properties for $\rho$ are identical, through the unitary change of frame $e^{-ih_\bn}$.\vspace{2mm}

\noindent \textbf{Lemma V.2:} \emph{The dynamics \eqref{eq:rhoKDecoherence} has an equilibrium of the form}
\begin{equation}\label{eq:eqform}
\bar{\rho}_\infty = \int_{-\alpha_c}^{\alpha_c} \mu(z) \, \q{z}\qd{z} \, dz
\end{equation}
\emph{where $\q{z}$ is a coherent state with $z \in \mathbb{R}$, $\alpha_c = 2\beta/(\kappa+\kappa_c)$, and}
\begin{equation}\label{eq:mu}
\mu(z) = \mu_0\, \frac{\left((\alpha_c^2 - z^2)^{\alpha_c^2}~e^{z^2}\right)^{(2 \kappa_c)/(\kappa+\kappa_c)}}{\alpha_c - z}\ ,
\end{equation}
\emph{with $\mu_0$ a normalization constant ensuring $\int_{-\alpha_c}^{\alpha_c} \mu(z) \, dz=1$.}

\underline{Proof:} One easily checks that \eqref{eq:rhoKDecoherence} leaves invariant the set of states of the form \eqref{eq:eqform} with arbitrary $\mu(z)$. Plugging this form into  \eqref{eq:rhoKDecoherence} and looking for a stationary solution $\mu(z)$, a projection on all real coherent states followed by a resolution by parts of the resulting ordinary differential equation  (with nonloncal terms!) yields the result; see Appendix B of \cite{CatsB} for details. \hfill $\square$\vspace{2mm}

In any case, $\mu(-\alpha_c) =0$. At the limit $\kappa_c \rightarrow 0$, the distribution $\mu(z)$ converges to a Dirac distribution at $z=\alpha_c$, implying:\vspace{2mm}

\noindent \textbf{Corollary V.3:} \emph{The dynamics \eqref{eq:rhoKDecoherence} without disturbance, i.e.~with $\kappa_c=0$, has an equilibrium of the form $\bar{\rho}_\infty = \q{\alpha_c}\qd{\alpha_c}$, for which $\rho = \rho_\infty$ is a two-component cat state of the form \eqref{def:cat}, more precisely:}
$$\rho_\infty = (\q{\tilde{\alpha_c}} + i \q{\minou\tilde{\alpha_c}}) \, (\q{\tilde{\alpha_c}} + i \q{\minou\tilde{\alpha_c}})\daag$$
\emph{with $\vert  \tilde{\alpha_c} \vert = \vert \alpha_c \vert$.}
\newline (This result is more easily obtained from \eqref{eq:cont-approx} directly.) \hfill $\square$\vspace{2mm}

For small $\kappa_c$, we have $\lim_{z\rightarrow \alpha_c}\mu(z)=+\infty$ so $\bar{\rho}_\infty$ is close to the coherent state $\q{\alpha_c}\qd{\alpha_c}$. We thus approach our goal. However, it remains to assess how/whether an arbitrary initial state \emph{converges towards} $\bar{\rho}_\infty$. For this we use the Corollary IV.3 developed above.\\

For starters, note that $\bar{\rho}_\infty$ has full rank as long as $\kappa_c \neq 0$. Indeed, the set of coherent states $\{ \q{\alpha} : \alpha \in [a,b] \subset \mathbb{R} \}$ makes an overcomplete basis of our Hilbert space, for any $b>a$. An annoying point is the infinite-dimensional setting. Since $n$ represents the energy of Fock state $\q{n}$ in units of $\hbar \omega$, physical arguments make it acceptable to limit ourselves to a finite-dimensional Hilbert space, spanned by $\q{0},\q{1},...,\q{\nmax}$: nobody physically expects the photon loss to transform the stable, coherent-state-stabilizing dynamics, into one where on average a significant fraction of the state is lost to infinitely high energies.
Within this approximation, we can then apply Corollary IV.3 directly to get the following.\vspace{2mm}

\noindent \textbf{Proposition V.4:} \emph{ Take $\nmax >0$ and denote by $\Pn$ the orthogonal projection onto $\mathcal{H}_{\nmax}= \text{span}\{\q{0},\q{1},...,\q{\nmax}\}$.  Set $\bnn=\Pn \bn \Pn$ and $\ban=\Pn \ba \Pn$ the truncation to $\mathcal{H}_{\bar n}$ of $\bn$ and $\ba$. The truncation of~\eqref{eq:rhoKDecoherence} to $\mathcal{H}_{\bar n}$ reads:
\begin{multline} \label{eq:rhoKDecoherenceTruncated}
	\tfrac{d}{dt} \rhon
 = \\
 \beta [ \ban\daag-\ban,\,  \rhon] - \tfrac{\kappa+\kappa_c}{2} (\bnn {\rhon} + {\rhon} \bnn - 2 \ban {\rhon} \ban\daag ) \\
 - \kappa_c (\ban {\rhon} \ban\daag - e^{i\pi\bnn}\ban {\rhon} \ban\daag e^{-i\pi\bnn} ) \; .
\end{multline}
where the support of $\rhon$ is included in $\mathcal{H}_{\bar n}$.
When \eqref{eq:rhoKDecoherenceTruncated} admits a full rank steady state $\rhonss$, the Lyapunov function $V_\text{Bures}(\rhon)$ applied to \eqref{eq:rhoKDecoherenceTruncated} is strict, showing that the truncated system globally converges towards $\rhonss$.}

\underline{Proof:} From Corollary IV.3, we have that  $V_\text{Bures}(\rhon)$ is a Lyapunov function for any Lindblad equation. It remains to examine the set where
\begin{multline*}
\tfrac{d}{dt} V_\text{Bures}= - \trace(\, [G,\ban]\, \rhonss \,[G,\ban]\daag \,) \\
  - \trace(\, [G,e^{i\pi\bnn}\ban]\, \rhonss\,[G,e^{i\pi\bnn}\ban]\daag \,)  \, = \, 0\, .
\end{multline*}
With $\rhonss$ positive definite, this occurs only when $G$ commutes with $\ban$ and with $e^{i\pi\bnn}\ban$; in fact one of them will be sufficient for our proof. Since $G$ is self-adjoint, we have the conditions $G \ban = \ban G$ and also $G \ban\daag = \ban\daag G$, which together yield $\bnn G = G \bnn$ since $\bnn=\ban\daag\ban$. Thus $G$ is diagonal in  basis $\{\q{0},\q{1},...,\q{\nmax}\}$, and as it commutes with $\ban$ it must in fact be a multiple of identity, $G= \lambda \bid_{\bar n}$.
Plugging this into \eqref{eq:LyapFunc1} and taking the trace shows that necessarily $\lambda=0$, thus $G=0$ which implies $\rhon = \rhonss$.
\hfill $\square$\vspace{2mm}

Thanks to the tools developed in Section IV, the proof above makes notably little use of system-specific properties. It can in fact be repeated verbatim to show the following.\vspace{2mm}

\noindent \textbf{Theorem V.5:} \emph{Consider any Lindblad-Kossakowski system \eqref{eq:Lindbladmap} which features one decoherence term $L_k$ proportional to any finite dimensional truncation $\ban$ of $\ba$ (defined in Proposition V.4) and has a full rank equilibrium $\rho_\infty$. Then $V_\text{Bures}(\rho)$ defined in Corollary IV.3 is a strict Lyapunov function for this system, so the latter globally converges to $\rho_\infty$.}\vspace{2mm}

Note that to have the same result with the DFS criterion~\cite{Davies1970,Spohn1977,Frigerio1978}, it is necessary to assume in addition that the subspace generated by all $L_k$ is Hermitian. This would be obtained e.g.~with a thermal bath at nonzero temperature. It seems worthwhile that Theorem V.5 gives a formal result also for environments at zero temperature, which is the ideal situation for physical experiments.


\section{Conclusion}\label{sec:conc}

On the basis of a set of contractive metrics for Kraus maps, we have systematically built a set of Lyapunov functions which decrease under \emph{any} Lindblad-Kossakowski evolution. In particular, unlike for many existing and popular convergence criteria such as fidelity, we obtain practical expressions for the rate of decay, allowing an efficient global investigation of mixed stationary states of the evolution. A major restriction of our approach is that it assumes an equilibrium of full rank. At the opposite side, pure state equilibria can usually be treated reasonably with fidelity-like indicators. An open question is thus how to efficiently treat systems with steady states of partial rank. We have illustrated our tool by proving convergence of an engineered reservoir scheme that stabilizes `Schr\"odinger cat' states of an electromagnetic field mode. For this we have introduced an approximation of finite dimensional Hilbert space. Infinite-dimensional Hilbert spaces, related to nonlocal PDE formulations, are typical of quantum systems like the basic harmonic oscillator. Systematic tools to alleviate the problems related to infinite dimension would thus be a welcome completion of this work.\\



\end{document}